# Product Calculus and Stokes' Theorem

M.G. Naber
mnaber@monroeccc.edu
Department of Science and Mathematics
Monroe County Community College
1555 S. Raisinville Rd
Monroe, Michigan, 48161-9746
November 30, 2023

**ABSTRACT**

In this paper the analogy between differential forms arising from integrals in additive calculus and *forms* arising from the integrals in product calculus is investigated. It is found that with an appropriate definition of scalar multiplication and vector addition a set of vector spaces can be constructed analogous to what is found in exterior calculus. A product differential is defined which allows for the product derivative version of closed and exact forms. The product differential also allows for a product integral version of Stokes' theorem (for scalar functions).

## I. INTRODUCTION

Within additive calculus, exterior calculus emerges, as Flanders puts it (pg. 1 Ref. 1), "forms…are the things which occur under integral signs". Differential forms provide an elegant means of expressing many results from differential geometry[1,2,3], the various theorems concerning multiple integrals (i.e. Green's, Gauss', Stokes', and the generalized Stokes' theorem[1,2]) and Maxwell's equations[1]. The question arises; what sort of form structure can be found in product integrals? In this paper the form structure implied by product integrals is examined. It is found to have vector spaces and a product differential analogous to what is found in exterior calculus. One might think of this as an exponential of an exterior algebra. The product form structure and associated differential are also used to produce a version of Stokes' theorem for product integrals.

In section 2 a brief review of product calculus will be given to fix notation and provide a convenient reference. This is by no means complete but is sufficient for the purposes of this paper. For a more complete treatment of product calculus the reader is encouraged to consult Refs. 4, 5, and 6. The notation and terminology for exterior calculus are from Refs. 1 and 2. In section 3 product forms are defined along with a product differential. With an appropriate definition of vector addition and scalar multiplication product forms are found to form a vector space. A vector product is defined which then allows for other vector spaces, analogous to the usual exterior product for differential forms. In section 4, after a brief review of simplices and chains, a version of Stokes' theorem for product integrals is given and shown to be true. In the appendix it is demonstrated how to compute product integrals for functions that take on negative values. The Einstein summation convention will be used,

$$\sum_i \frac{\partial y^j}{\partial x^i} dx^i = \frac{\partial y^j}{\partial x^i} dx^i.$$

## II. PRODUCT CALCULUS

Calculus can be thought of as a branch of mathematics where operators are constructed out of sums, differences, products, and quotients of function values infinitesimally far apart. The two most commonly used are the additive derivative operator (slope) and the additive integration operator (area under a curve). Rather than considering additive changes in a function, multiplicative changes can be considered. The multiplicative derivative (also called a product derivative in the literature, c.f. Refs. 4, 5, and 6) is defined as,

$$q_x f(x) = \lim_{\delta \to 0} \left( \frac{f(x+\delta)}{f(x)} \right)^{1/\delta} = e^{\frac{\partial_x f(x)}{f(x)}}. \tag{1}$$

$q_x$ will be the notation for the multiplicative derivative just as $\partial_x$ is the notation for the additive derivative. Just as the additive derivative provides information about the additive change in a function for points $x$ near $c$, the multiplicative derivative provides information about the multiplicative change in a function near $c$. For $x$ very near $c$,

$$f(x) \approx f(c)(q_x f(x)|_{x=c})^{(x-c)}. \tag{2}$$

Multiplicative derivatives have found use in a variety of areas. A few examples are, physics and approximation schemes[5], finance and economics[4,7], and image processing[8].

Similarly, to the additive integral, a multiplicative integral can be defined (c.f. Sec. 2.6 of Ref. 4). Let $f$ be defined and positive on $[a, b]$ (See the appendix for the case of $f(x) < 0$) and consider a partition, $a = x_0 < x_1 < \cdots < x_{n-1} < x_n = b$. Define $\Delta_k = x_k - x_{k-1}$, let $\Delta = max\{\Delta_k, k = 1 \ldots n\}$, and let $c_k$ be a number such that $x_{k-1} < c_k < x_k$, then,

$$\prod_a^b f(x)^{dx} = \lim_{\Delta \to 0} \prod_{k=1}^n f(c_k)^{\Delta_k} = e^{\int_a^b \ln(f(x)) dx}. \tag{3}$$

The first part of the above equation is notation, the upper-case pi replaces the elongated S used for sum integration. The middle part is how it might be computed numerically, or how we might draw a picture to see what is happening. The third part is how it is related to an additive integral. Note that if $f(x)$ is matrix valued the last part of Eq. (9) is not valid and the order in which the product is taken within the limit must also be specified (c.f. chapter 1 of Ref. 6). This product integral has found uses in a wide variety of areas. Some examples are, biomedical imaging[9], probability and statistics[10], economics and chaos[11], and, if $f(x)$ is allowed to take on matrix values, the study of differential equations[6] and path integrals and quantum mechanics[12].

Notice that the product derivative is the left inverse for the above defined product integral,

$$q_x\left(\prod_a^x f(t)^{dt}\right) = f(x), \tag{4}$$

$$\prod_a^x (q_t f(t))^{dt} = \frac{f(x)}{f(a)}, \tag{5}$$

$$\prod (q_x f(x))^{dx} = f(x)c. \tag{6}$$

The above can all be viewed as statements of the fundamental theorem of Multiplicative calculus (c.f. pg. 20 of Ref. 6). These would be the multiplicative analogues of,

$$\partial_x \int_a^x f(t)dt = f(x), \tag{7}$$

$$\int_a^x \partial_t f(t)dt = f(x) - f(a), \tag{8}$$

$$\int \partial_x f(x)dx = f(x) + c. \tag{9}$$

The other commonly used product integral is the Volterra product integral[13,14]. In a formal manner, this can be obtained from Eq. (3). Express the integrand as an exponential,

$$\prod_a^b f(x)^{dx} = \lim_{\Delta \to 0} \prod_{k=1}^n e^{\ln(f(c_k))\Delta_k}. \tag{10}$$

Now expand the exponential,

$$= \lim_{\Delta \to 0} \prod_{k=1}^n \left(1 + \ln(f(c_k))\Delta_k + \ln(f(c_k))^2 \frac{(\Delta_k)^2}{2} + \cdots\right). \tag{11}$$

If the square and higher powers of the infinitesimal, $\Delta_k$, are taken to be vanishingly small,

$$= \lim_{\Delta \to 0} \prod_{k=1}^n (1 + \ln(f(c_k))\Delta_k). \tag{12}$$

Now replace $\ln(f(x))$ with $g(x)$, and the Volterra product integral is obtained,

$$\prod_a^b (1 + g(x)dx) = \lim_{\Delta \to 0} \prod_{k=1}^n (1 + g(c_k)\Delta_k) = e^{\int_a^b g(x)dx}. \tag{13}$$

The left inverse for the Volterra integral is the so-called logarithmic derivative,

$$\ell_x f(x) = \frac{\partial_x f(x)}{f(x)}. \tag{14}$$

The Volterra analogues to Eqs. (4), (5), (6), and (7), (8), (9), are then,

$$\ell_x \prod_a^x (1 + g(t)dt) = g(t), \tag{15}$$

$$\prod_a^x (1 + \ell_t g(t)dt) = \frac{g(x)}{g(a)}, \tag{16}$$

$$\prod (1 + \ell_x g(x)dx) = g(x)c. \tag{17}$$

## III. PRODUCT FORMS

To motivate the analogy with differential forms, consider some integrals in $\mathbb{R}^3$. Recall that an additive integral along a curve generates a 1-form, $\alpha$,

$$\int_{curve} A dx + B dy + C dz, \tag{18}$$

$$\alpha = A dx + B dy + C dz. \tag{19}$$

Similarly, a product integral along a curve would generate a 1-product form, $\lambda$,

$$\prod_{curve} A^{dx} B^{dy} C^{dz}, \tag{20}$$

$$\lambda = (A)^{dx}(B)^{dy}(C)^{dz}. \tag{21}$$

Carrying on, 2 and 3-product forms, $\mu$ and $\nu$, arise from surface and volume product integrals,

$$\prod_{surface} D^{dxdy} E^{dydz} F^{dzdx}, \tag{22}$$

$$\mu = (D)^{dxdy}(E)^{dydz}(F)^{dzdx}, \tag{23}$$

$$\prod_{volume} G^{dxdydz}, \tag{24}$$

$$\nu = (G)^{dxdydz}. \tag{25}$$

The product form version of the exterior derivative, $d$, is denoted by $q$, and is defined as,

$$qf = (q_x f)^{dx}(q_y f)^{dy}(q_z f)^{dz}, \tag{26}$$

Or, in $\mathbb{R}^n$,

$$qf = \prod_{i=1}^{n}(q_{x_i} f)^{dx_i} = e^{d \ln(f)}. \tag{27}$$

The product form version of the Poincaré lemma is then,

$$q^2 f = e^{d^2 \ln(f)} = e^0 = 1. \tag{28}$$

Notice that for product forms, 1 will play the same role as 0 does for differential forms.

As an example, let $\lambda$ be a 1-product form in $\mathbb{R}^3$, $\lambda = (A)^{dx}(B)^{dy}(C)^{dz}$. The $q$ differential of $\lambda$ is then,

$$q\lambda = \left(e^{\partial_x \ln(C) - \partial_z \ln(A)}\right)^{dzdx} \left(e^{\partial_y \ln(A) - \partial_x \ln(B)}\right)^{dxdy} \left(e^{\partial_z \ln(B) - \partial_y \ln(C)}\right)^{dydz}. \tag{29}$$

The second $q$ differential is,

$$q^2\lambda = \left(e^{\partial_y(\partial_x \ln(C) - \partial_z \ln(A)) + \partial_z(\partial_y \ln(A) - \partial_x \ln(B)) + \partial_x(\partial_z \ln(B) - \partial_y \ln(C))}\right)^{dxdydz} = 1. \tag{30}$$

Let $\mu$ be a 2-product form in $\mathbb{R}^3$, $\mu = (A)^{dxdy}(B)^{dydz}(C)^{dzdx}$. The $q$ differential is then,

$$q\mu = \left(e^{\partial_x \ln(B) + \partial_y \ln(C) + \partial_z \ln(A)}\right)^{dxdydz}, \tag{31}$$

and,

$$q^2\mu = 1. \tag{32}$$

A product form $\mu$ will be called closed if $q\mu = 1$. A product form, $\omega$, of order $p$, will be called exact if there is a product form of order $p - 1$, $\eta$, such that $\omega = q\eta$.

In a more formal setting, product forms are defined below. Let $\mathbb{R}$ be the field of real numbers, $i, j, k, \ldots$ be indices $\{1, 2, \ldots, n\}$. Let $a_i, b_i, \ldots$ be in $\mathbb{R}^+$ (i.e. real numbers greater than zero). Let $\alpha_1, \alpha_2, \ldots, \alpha_n$ be basis vectors for an $n$ dimensional exterior algebra with the usual rules of vector addition and scalar multiplication (see chapter 2 of Ref. 1). Then, $\mathbb{P}$ is an $n$-dimensional vector space defined below. $\mathbb{P}^p$ denotes a related vector space of dimension $\binom{n}{p}$ for $p = 1, 2, \ldots, n$ ($\mathbb{P}^1 = \mathbb{P}$). Let $\alpha, \beta,$ and $\gamma$ be vectors in $\mathbb{P}$,

$$\alpha = \prod_{i=1}^{n}(a_i)^{\alpha_i}, \quad \beta = \prod_{i=1}^{n}(b_i)^{\alpha_i}, \quad \gamma = \prod_{i=1}^{n}(c_i)^{\alpha_i}. \tag{33}$$

For this set define vector addition to be,

$$\alpha \oplus \beta = \left(\prod_{j=1}^{n}(a_j)^{\alpha_j}\right) \oplus \left(\prod_{k=1}^{n}(b_k)^{\alpha_k}\right) = \prod_{i=1}^{n}(a_i b_i)^{\alpha_i}. \tag{34}$$

The identity element for vector addition is denoted by $I$ and is defined as,

$$I = \prod_{i=1}^{n}(1)^{\alpha_i}. \tag{35}$$

Scalar multiplication is defined to be ($a = 1$ is the identity),

$$a \odot \alpha = a \odot \prod_{i=1}^{n}(a_i)^{\alpha_i} = \prod_{i=1}^{n}(a_i{}^a)^{\alpha_i}. \tag{36}$$

The inverse for a product form $\alpha$ is then,

$$inv(\alpha) = \prod_{i=1}^{n}\left(\frac{1}{a_i}\right)^{\alpha_i} = (-1) \odot \prod_{i=1}^{n}(a_i)^{\alpha_i}. \tag{37}$$

The above defined vector addition, Eq. (34), and scalar multiplication, Eq. (36), have the appropriate associative, commutative, and distributive properties for $\mathbb{P}$ to form a vector space (pg. 43 Ref. 2). The other $\mathbb{P}^p$ also form vector spaces with elements of the form,

$$\mathbb{P}^2 \rightarrow \prod_{i=1}^{n-1}\left(\prod_{j=1+i}^{n}(a_{ij})^{\alpha_i \wedge \alpha_j}\right), \tag{38}$$

$$\mathbb{P}^3 \rightarrow \prod_{i=1}^{n-2}\left(\prod_{j=1+i}^{n-1}\left(\prod_{k=1+j}^{n}(a_{ijk})^{\alpha_i \wedge \alpha_j \wedge \alpha_k}\right)\right), \tag{39}$$

$$\vdots$$

$$\mathbb{P}^n \rightarrow (a)^{\alpha_1 \wedge \alpha_2 \wedge \ldots \wedge \alpha_n}. \tag{40}$$

Objects in $\mathbb{P}^2$ can be constructed from objects in $\mathbb{P}^1$ using the following definition of a vector product (denoted by $\wedge_p$) for product forms. Let $(a)^{\alpha_1}$ and $(b)^{\alpha_2}$ be 1-product forms, then,

$$(a)^{\alpha_1} \wedge_p (b)^{\alpha_2} = (ab)^{\alpha_1 \wedge \alpha_2}. \tag{41}$$

Notice,

$$(a)^{\alpha_i} \wedge_p (b)^{\alpha_j} = (-1) \odot \left((b)^{\alpha_j} \wedge_p (a)^{\alpha_i}\right). \tag{42}$$

In general, if $\alpha$ is a $p$-product form and $\beta$ is a $q$-product form then,

$$\alpha \wedge_p \beta = (-1)^{pq} \odot (\beta \wedge_p \alpha). \tag{43}$$

This vector product is also associative and distributive.

$$\alpha \wedge_p (\beta \wedge_p \gamma) = (\alpha \wedge_p \beta) \wedge_p \gamma \tag{44}$$

$$\alpha \wedge_p (\beta \oplus \gamma) = (\alpha \wedge_p \beta) \oplus (\alpha \wedge_p \gamma) \tag{45}$$

As an example, let $\alpha$, and $\beta$, be 1-product forms as given in Eq. (33). Then,

$$\alpha \wedge_p \beta = \left(\prod_{i=1}^n (a_i)^{\alpha_i}\right) \wedge_p \left(\prod_{j=1}^n (b_j)^{\alpha_j}\right) \tag{46}$$

$$= \left((a_1)^{\alpha_1} \wedge_p ((b_1)^{\alpha_1}(b_2)^{\alpha_2} \ldots (b_n)^{\alpha_n})\right)\left((a_2)^{\alpha_2} \wedge_p ((b_1)^{\alpha_1}(b_2)^{\alpha_2} \ldots (b_n)^{\alpha_n})\right)$$
$$\ldots \left((a_n)^{\alpha_n} \wedge_p ((b_1)^{\alpha_1}(b_2)^{\alpha_2} \ldots (b_n)^{\alpha_n})\right) \tag{47}$$

Noting that, $\alpha_k \wedge \alpha_k = 0$, and gathering like terms, Eq. (47) can be reduced to,

$$\left(\prod_{i=2}^n \left(\frac{a_1 b_i}{a_i b_1}\right)^{\alpha_1 \wedge \alpha_i}\right)\left(\prod_{i=3}^n \left(\frac{a_2 b_i}{a_i b_2}\right)^{\alpha_2 \wedge \alpha_i}\right) \ldots \left(\left(\frac{a_{n-1} b_n}{a_n b_{n-1}}\right)^{\alpha_{n-1} \wedge \alpha_i}\right) \tag{48}$$

The $q$ differential acts with the vector product in much the same way as the $d$ differential acts with the usual exterior product. Here, let $\alpha$ be an $n$-product form, then,

$$q(\alpha \wedge_p \beta) = (q\alpha) \wedge_p \beta \oplus (-1)^n \odot \left(\alpha \wedge_p (q\beta)\right). \tag{49}$$

**IV. STOKES' THEOREM**

For additive integrals, Stokes' theorem can be stated as (c.f. pg. 64 of Ref. 1). Let $\omega$ be a $p$-form on a manifold $M$ and $c$ a $(p+1)$-chain. Then,

$$\int_{\partial c} \omega = \int_{c} d\omega. \tag{50}$$

The analogous statement for product integrals is,

$$\prod_{\partial c} \alpha = \prod_{c} q\alpha, \tag{51}$$

where $\alpha$ is a $p$-product form. The proof of Stokes' theorem for differential forms within additive integrals relies on the use of simplices, chains, and the transformation properties of differential forms (c.f. pgs. 64-66 of Ref. 1). The proof of Stokes' theorem for product integrals follows essentially the same steps. The same notation as Flanders[1] shall be used for simplices and chains.

Recall that an $n$-simplex is a closed convex hull of $n+1$ independent points taken in a specific order, $\boldsymbol{s} = (P_0, P_1, \ldots, P_n)$. This set consists of all points, $P$, such that,

$$P = t_0 P_0 + t_1 P_1 + \cdots + t_n P_n, \quad t_i \geq 0, \quad \sum_{i=0}^{n} t_i = 1. \tag{52}$$

The boundary of a simplex is the signed sum of all possible simplices one dimension lower.

$$\partial \boldsymbol{s} = \partial(P_0, P_1, \ldots, P_n) = \sum_{i=0}^{n} (-1)^i (P_0, P_1, \ldots \hat{P}_i, \ldots, P_n) \tag{53}$$

Where $\hat{P}_i$ indicates that $P_i$ does not appear in the term. The standard points in $\mathbb{R}^n$ to construct the standard $n$-simplex are,

$$\begin{aligned} R_0 &= (0, \ldots, 0), \\ R_1 &= (1, 0, \ldots, 0), \\ &\vdots \\ R_n &= (0, \ldots, 0, 1). \end{aligned} \tag{54}$$

An $n$-chain is a sum of simplices (the $a_i$ are constants),

$$c = \sum a_i \boldsymbol{s}_i. \tag{55}$$

The other tool that is needed is the behavior of forms under a mapping. For differential forms the mapping is constructed as follows (see sec. 3.3 of Ref. 1). Let $\boldsymbol{M}$ be an $m$-dimensional manifold, $\boldsymbol{U}$ a domain in $\boldsymbol{M}$ with coordinates denoted by $x^1, \ldots, x^m$. Let $\boldsymbol{N}$ be an $n$-dimensional manifold, $\boldsymbol{V}$ domain in $\boldsymbol{N}$ with coordinates denoted by $y^1, \ldots, y^n$, and $\phi$ a smooth mapping that takes $U$ into $V$. The mapping $\phi$ induces a mapping $\phi_k^*$ that maps $k$-forms on $V$ into $k$-forms on $U$. For 1-forms the mapping is given by,

$$\phi_1^*\left(\sum_i a_i(y)dy^i\right) = \sum_i a_i(y(x))\frac{\partial y^i}{\partial x^j}dx^j. \tag{56}$$

Note that $\phi_k^*$ has the following properties. Let $\alpha$ and $\beta$ be $m$-forms and $\gamma$ be an $n$-form on $V$.

$$\phi_m^*(\alpha + \beta) = \phi_m^*(\alpha) + \phi_m^*(\beta) \tag{57}$$

$$\phi_{m+n}^*(\alpha \wedge \gamma) = \phi_m^*(\alpha) \wedge \phi_n^*(\gamma) \tag{58}$$

$$\phi_{m+1}^*(d\alpha) = d\phi_m^*(\alpha) \tag{59}$$

For product forms the mapping is very similar and will be denoted by $\phi_k^{*p}$. For 1-product forms the mapping is given by,

$$\phi_1^{*p}\left(\prod_i (a_i)^{dy^i}\right) = \prod_i (a_i)^{\frac{\partial y^i}{\partial x^j}dx^j} = \prod_i \left(a_i^{\frac{\partial y^i}{\partial x^j}}\right)^{dx^j}. \tag{60}$$

For 2-product forms the mapping is given by,

$$\phi_2^{*p}\left(\prod_{i,j} (a_{ij})^{dy^i dy^j}\right) = \prod_{i,j} \left(a_{ij}^{\frac{\partial y^i \partial y^j}{\partial x^k \partial x^l}}\right)^{dx^j dx^l}. \tag{61}$$

The product form formula analogues to Eqs. (57), (58), and (59) are then,

$$\phi_k^{*p}(\alpha \oplus \beta) = \phi_k^{*p}(\alpha) \oplus \phi_k^{*p}(\beta), \tag{62}$$

$$\phi_{k+l}^{*p}(\alpha \wedge_p \gamma) = \phi_k^{*p}(\alpha) \wedge_p \phi_l^{*p}(\gamma), \tag{63}$$

$$q(\phi_k^{*p}\alpha) = \phi_{k+1}^{*p}(q\alpha). \tag{64}$$

Stokes' theorem for product forms can now be stated and proved. Let $M$ be a $m$-dimensional manifold. Let $c$ be a $(n+1)$-chain, $c = \sum a_i \sigma_i$ where the $a_i$ are real constants and the $\sigma_i$ $(n+1)$-simplices. Let $\alpha$ be a $n$-product form on $M$. Then,

$$\prod_{\partial c} \alpha = \prod_c q\alpha, \tag{65}$$

or,

$$\prod_{\Sigma a_i(\partial \sigma_i)} \alpha = \prod_{\Sigma a_i \sigma_i} q\alpha. \tag{66}$$

As with additive integrals, $c$ is a sum of $(n+1)$-simplices, hence, it is sufficient to prove only for a single simplex, i.e.

$$\prod_{\partial \sigma} \alpha = \prod_{\sigma} q\alpha, \tag{67}$$

where $\sigma$ is an $(n+1)$-simplex. Let $\phi$ be a mapping from $W$, an open neighborhood of $\bar{s}^{n+1}$ in $\mathbb{R}^{n+1}$ to $M$, such that $\phi(\bar{s}^{n+1}) = \sigma$, and $\bar{s}^{n+1}$ is the standard $(n+1)$-simplex (recall, $\bar{s}^{n+1}$ is the set of all points such that $x^i \geq 0$ and $\sum_{i=1}^{n+1} x^i \leq 1$). Using this representation, Eq. (67) is reduced to a product integral in Euclidean space,

$$\prod_{\partial \bar{s}^{n+1}} \phi_n^{*p}(\alpha) = \prod_{\bar{s}^{n+1}} \phi_{n+1}^{*p}(q\alpha). \tag{68}$$

Where the induced mapping, $\phi_n^{*p}$, is defined above and commutes with the $q$ differential, according to Eq. (64), to give,

$$\prod_{\partial \bar{s}^{n+1}} \phi_n^{*p}(\alpha) = \prod_{\bar{s}^{n+1}} q\left(\phi_n^{*p}(\alpha)\right). \tag{69}$$

The product form, $\phi_n^{*p}(\alpha)$, in the integrals above is a product form in $\mathbb{R}^n$. Write this product form as

$$\phi_n^{*p}(\alpha) = \eta. \tag{70}$$

The proof is then reduced to showing,

$$\prod_{\partial \bar{s}^{n+1}} \eta = \prod_{\bar{s}^{n+1}} q\eta, \tag{71}$$

where,

$$\eta = \prod_{i=1}^{n}(A_i)^{dx^1\ldots d\hat{x}^i\ldots dx^{n+1}}, \tag{72}$$

($d\hat{x}^i$ indicates the $dx^i$ term that is missing). As in the additive integral version of Stokes' theorem it is sufficient to consider the case (where a convenient choice of the order of the coordinate labels has been made, and $A(x) = A(x^1, x^2, \ldots, x^n)$),

$$\eta = (A(x))^{dx^1\ldots dx^n}, \tag{73}$$

and then,

$$q\eta = (-1)^n \odot \left(q_{x^{n+1}} A(x)\right)^{dx^1 \ldots dx^{n+1}}. \tag{74}$$

Substituting Eq. (74) into the right-hand side of Eq. (71) gives,

$$\prod_{\bar{s}^{n+1}} q\eta = \prod_{\bar{s}^{n+1}} (-1)^n \odot \left(q_{x^{n+1}} A(x)\right)^{dx^1 \ldots dx^{n+1}}, \tag{75}$$

$$= (-1)^n \odot \prod_{\bar{s}^{n+1}} \left(q_{x^{n+1}} A(x)\right)^{dx^1 \ldots dx^{n+1}}. \tag{76}$$

Now break the product integral into two parts, the first part will cover the variables $x^1 \ldots x^n$ and the second part over $x^{n+1}$ (see the comments above Eq. (68)),

$$= (-1)^n \odot \prod_{\substack{x^i \geq 0 \\ \sum_{i=1}^n x^i \leq 1}} \left( \prod_{x^{n+1}=0}^{1-\sum_{i=1}^n x^i} \left(q_{x^{n+1}} A(x)\right)^{dx^{n+1}} \right)^{dx^1 \ldots dx^n}, \tag{77}$$

Using Eq. (5), Eq. (77) reduces to,

$$= (-1)^n \odot \prod_{\substack{x^i \geq 0 \\ \sum_{i=1}^n x^i \leq 1}} \left( \frac{A\left(x^1, \ldots x^n, 1 - \sum_{i=1}^n x^i\right)}{A(x^1, \ldots x^n, 0)} \right)^{dx^1 \ldots dx^n}. \tag{78}$$

Next, investigate the integral on the left-hand side of Eq. (71). Recall,

$$\bar{s}^{n+1} = (R_0, R_1, \ldots, R_{n+1}). \tag{79}$$

The boundary is the signed sum of all the faces,

$$\partial \bar{s}^{n+1} = (R_1, \ldots, R_{n+1}) + (-1)^{n+1}(R_0, R_1, \ldots, R_n) + \cdots \tag{80}$$

On all the faces not explicitly listed, $\eta = 1$, as at least one of the $x^1, \ldots, x^n$ are constants.

$$\prod_{\partial \bar{s}^{n+1}} \eta = \left( \prod_{(R_1,\ldots,R_{n+1})} \eta \right) \oplus (-1)^{n+1} \odot \prod_{(R_0,\ldots,R_n)} \eta \tag{81}$$

The domain for the second integral on the right-hand side of Eq. (81) has $x^{n+1} = 0$,

$$\prod_{\partial \bar{s}^{n+1}} \eta = \left( \prod_{(R_1,\ldots,R_{n+1})} \eta \right) \oplus (-1)^{n+1} \odot \prod_{(R_0,\ldots,R_n)} \left(A(x^1, \ldots, x^n, 0)\right)^{dx^1 \ldots dx^n}. \tag{82}$$

In the first integral on the right-hand side of Eq. (82), $(R_1, \dots, R_{n+1})$ can be projected along $x^{n+1}$ to the $(R_1, \dots, R_n, R_0)$-plane. This gives,

$$\prod_{(R_1,\dots,R_{n+1})} \eta = \prod_{(R_1,\dots,R_n,R_0)} \left(A\left(x^1, \dots, x^n, 1 - \sum_{i=1}^{n} x^i\right)\right)^{dx^1 \dots dx^n}. \tag{83}$$

Or, if the points in the simplex are reordered,

$$= (-1)^n \odot \prod_{(R_0,\dots,R_n)} \left(A\left(x^1, \dots, x^n, 1 - \sum_{i=1}^{n} x^i\right)\right)^{dx^1 \dots dx^n}. \tag{84}$$

Combining what has been obtained for the two integrals on the right-hand side of Eq. (81) gives,

$$\prod_{\partial \bar{s}^{n+1}} \eta = (-1)^n \odot \prod_{(R_0,\dots,R_n)} \left(A\left(x^1, \dots, x^n, 1 - \sum_{i=1}^{n} x^i\right)\right)^{dx^1 \dots dx^n}$$
$$\oplus (-1)^{n+1} \odot \prod_{(R_0,\dots,R_n)} \left(A(x^1, \dots, x^n, 0)\right)^{dx^1 \dots dx^n}. \tag{85}$$

Factoring out the $(-1)^n$ and combining integrals gives,

$$\prod_{\partial \bar{s}^{n+1}} \eta = (-1)^n \odot \prod_{(R_0,\dots,R_n)} \left(\frac{A\left(x^1, \dots, x^n, 1 - \sum_{i=1}^{n} x^i\right)}{A(x^1, \dots, x^n, 0)}\right)^{dx^1 \dots dx^n} \tag{86}$$

Or,

$$= (-1)^n \odot \prod_{\substack{x^i \geq 0 \\ \sum_{i=1}^{n} x^i \leq 1}} \left(\frac{A\left(x^1, \dots x^n, 1 - \sum_{i=1}^{n} x^i\right)}{A(x^1, \dots x^n, 0)}\right)^{dx^1 \dots dx^n}. \tag{87}$$

Eq. (87) is seen to be the same as Eq. (78) which proves Stokes' theorem for product integrals.

## V. CONCLUSION

In this paper product forms were found by examining product integrals in analogy with exterior calculus. After a suitable form of vector addition and multiplication and scalar multiplication were defined product forms were found to form vector spaces and have analogous properties to the vector spaces associated with exterior calculus. Having identified this form structure, and a product differential, a product integral version of Stokes' theorem was shown to be true. Future studies would be concerned with adapting the formalism presented here to matrix valued functions which would have application to systems of differential equations. As side result, presented in the appendix, demonstrates how to compute product integrals for functions that take on negative values. The presented method has the properties one would expect.

**ACKNOWLEDGMENTS**

The author would like to thank J. P. Krisch for critically reading an early version of this paper and making helpful suggestions towards improving the exposition of this article.

**CONFLICTS OF INTEREST**

The author has no conflicts of interest to disclose.

**DECLARATIONS OF INTEREST**: none.

**DATA AVAILABILLITY**

Data sharing is not applicable to this article as no new data was created or analyzed in this study.

**APPENDIX: GEOMETRIC MEAN OF NEGATIVE NUMBERS**

Just as the additive integral can be used to compute the additive mean of a function over an interval, namely,

$$\frac{1}{b-a}\int_a^b f(t)dt. \tag{88}$$

The geometric integral can be used to compute the geometric mean of a function over an interval,

$$\left(\prod_a^b f(x)^{dx}\right)^{\frac{1}{b-a}}. \tag{89}$$

What if $f(x) < 0$ on the interval $[a, b]$? How is the product integral computed? Consider the following theorem: Let $f(x) < 0 \ \forall x \in [a, b]$ and let $|f(x)|$ be product integrable on $[a, b]$. Then the following is true,

$$\prod_a^b f(x)^{dx} = e^{i\pi(b-a)} \prod_a^b |f(x)|^{dx}. \tag{900}$$

Proof: Let $\{x_0, x_1, ..., x_n\}$ be a partition of $[a, b]$ such that $x_0 = a$, $x_n = b$, and $x_{k-1} < x_k$. Let $c_k$ be any point in the interval $[x_{k-1}, x_k]$ and let $\Delta_k = x_k - x_{k-1}$, be the width of the interval. Let $\Delta$ be the maximum value of the $\Delta_k$. Recall that $-1 = e^{i\pi}$, then, on $[a, b]$ $f(x)$ can be written as,

$$f(x) = e^{i\pi}|f(x)|. \tag{91}$$

The product integral can be written as

$$\prod_a^b f(x)^{dx} = \lim_{\Delta \to 0} \prod_{k=1}^n e^{i\pi \Delta_k}|f(c_k)|^{\Delta_k}. \tag{92}$$

The product can then be split in two,

$$\prod_a^b f(x)^{dx} = \lim_{\Delta \to 0} \prod_{k=1}^n e^{i\pi \Delta_k} \prod_{k=1}^n |f(c_k)|^{\Delta_k}, \tag{93}$$

or,

$$\prod_a^b f(x)^{dx} = \left( \lim_{\Delta \to 0} \prod_{k=1}^n e^{i\pi \Delta_k} \right) \left( \lim_{\Delta \to 0} \prod_{k=1}^n |f(c_k)|^{\Delta_k} \right). \tag{94}$$

The second product is just the product integral of $|f(x)|$.

$$\prod_a^b f(x)^{dx} = \left( \lim_{\Delta \to 0} \prod_{k=1}^n e^{i\pi \Delta_k} \right) \prod_a^b |f(x)|^{dx} \tag{95}$$

Now consider the remaining limit by itself,

$$\lim_{\Delta \to 0} \prod_{k=1}^n e^{i\pi \Delta_k} = \lim_{\Delta \to 0} e^{i\pi \sum_{k=1}^n \Delta_k} = e^{i\pi(b-a)}. \tag{96}$$

Hence,

$$\prod_a^b f(x)^{dx} = e^{i\pi(b-a)} \prod_a^b |f(x)|^{dx}. \tag{97}$$

Notice that the formula for product integration for functions that are negative on an interval gives the expected result. Let $f(x) < 0$ for all $x$ in the interval $[a, b]$, then the geometric mean value of the function is,

$$\left( \prod_a^b f(x)^{dx} \right)^{\frac{1}{b-a}} = -\left( \prod_a^b |f(x)|^{dx} \right)^{\frac{1}{b-a}}. \tag{98}$$

As an example, compute the geometric mean for $sin(x)$ over the interval $[0, 2\pi]$.

$$\left(\prod_0^{2\pi} sin(x)^{dx}\right)^{\frac{1}{2\pi}} = \left(\left(\lim_{a\to 0^+}\prod_a^{\pi-a} sin(x)^{dx}\right)\left(\lim_{a\to 0^+} e^{i\pi(\pi-2a)}\prod_{\pi+a}^{2\pi-a}|sin(x)|^{dx}\right)\right)^{\frac{1}{2\pi}} \quad (99)$$

$$= i\left(\left(\lim_{a\to 0^+}\prod_a^{\pi-a} sin(x)^{dx}\right)\left(\lim_{a\to 0^+}\prod_{\pi+a}^{2\pi-a}|sin(x)|^{dx}\right)\right)^{\frac{1}{2\pi}}, \quad (100)$$

$$= i\left(\left(\lim_{a\to 0^+} e^{\int_a^{\pi-a} ln(sin(x))dx}\right)\left(\lim_{a\to 0^+} e^{\int_{\pi+a}^{2\pi-a} ln(|sin(x)|)dx}\right)\right)^{\frac{1}{2\pi}}, \quad (101)$$

$$= i\left(\left(e^{-\pi ln(2)}\right)\left(e^{-\pi ln(2)}\right)\right)^{\frac{1}{2\pi}}, \quad (102)$$

$$= \frac{i}{2}. \quad (103)$$

The geometric mean of $sin(x)$ on the interval $(0, \pi)$ is $1/2$ and likewise for the geometric mean of $|sin(x)|$ on the interval $(\pi, 2\pi)$. Over the full interval $(0, 2\pi)$, excluding $x = \pi$, sine is positive on half the interval, and it is negative on the other half, so the purely imaginary result makes sense.